\documentstyle[amscd,amssymb,verbatim,11pt]{amsart}

\theoremstyle{plain}
\newtheorem{Prop}{Proposition}[section]
\newtheorem{Thm}[Prop]{Theorem}
\newtheorem{Cor}[Prop]{Corollary}

\newtheorem{Lem}[Prop]{Lemma}

\theoremstyle{definition}
\newtheorem{Def}[Prop]{Definition}

\theoremstyle{remark}
\newtheorem{Rem}[Prop]{Remark}
\newtheorem{ex}[Prop]{Example}
\newtheorem{Problem}[Prop]{\bf Problem}

\def\dim{\mathop{\roman{dim}}}
\def\uudim{\dim_u}
\def\int{\mathop{\roman{int}}}

\def\1{^{-1}}

\def\Cone{\text{Cone}}

\def\dim{\text{dim}}
\def\diam{\text{diam}}
\def\id{\text{id}}
\def\Lip{\text{Lip}}

\def\dist{\text{dist}}

\def\uD{{\mathcal D}}

\def\LL{{L}}

\def\MM{{\mathcal M}}
\def\NN{{\mathbb N}}
\def\RR{{\mathbb R}}

\def\ZZ{{\mathbb Z}}

\errorcontextlines=0 \numberwithin{equation}{section}

\begin{document}
\title[Dimension zero at all scales]{Dimension zero at all scales}

\author{N.~Brodskiy}
\address{University of Tennessee, Knoxville, TN 37996, USA}
\email{brodskiy@@math.utk.edu}

\author{J.~Dydak}
\address{University of Tennessee, Knoxville, TN 37996, USA}
\email{dydak@@math.utk.edu}

\author{J.~Higes}
\address{Departamento de Geometr\'{\i}a y Topolog\'{\i}a,
Facultad de CC.Matem\'aticas. Universidad Complutense de Madrid.
Madrid, 28040 Spain} \email{josemhiges@@yahoo.es}

\author{A.~Mitra}
\address{University of Tennessee, Knoxville, TN 37996, USA}
\email{ajmitra@@math.utk.edu}




\begin{abstract}
We consider the notion of dimension in four categories: the
category of (unbounded) separable metric spaces and (metrically
proper) Lipschitz maps, and the category of (unbounded) separable
metric spaces and (metrically proper) uniform maps. A unified
treatment is given to the large scale dimension and the small
scale dimension. We show that in all categories a space has
dimension zero if and only if it is equivalent to an ultrametric
space. Also, 0-dimensional spaces are characterized by means of
retractions to subspaces. There is a universal zero-dimensional
space in all categories. In the Lipschitz Category spaces of
dimension zero are characterized by means of extensions of maps to
the unit 0-sphere. Any countable group of asymptotic dimension
zero is coarsely equivalent to a direct sum of cyclic groups. We
construct uncountably many examples of coarsely inequivalent
ultrametric spaces.
\end{abstract}

\maketitle

\tableofcontents

\section{Introduction}

Asymptotic dimension is one of the most important asymptotic
invariants of metric spaces introduced by Gromov~\cite{Gro asym
invar}. There are several notions of large scale dimension
introduced later~\cite{Dran asym top, Dran asind, Brod-Dydak}. The
asymptotic dimension of Gromov is known to be the largest and in
case it is finite all dimensions coincide. These dimensions also
coincide when one of them is zero, but it is still unknown if an
example of space exists with one of these dimensions finite but
the asymptotic dimension of Gromov infinite. The notion of
asymptotic dimension can be introduced for any set with coarse
structure~\cite{Roe lectures} (or a ballean~\cite{ProtasovSurvey,
Banakh Protasov}) but in this paper we consider separable metric
spaces only.

Our attempts to find the small scale analogs of large scale
dimensions brought us to an idea of macroscopic and microscopic
functors on a category of metric spaces: given a metric space
$(X,d)$ and $\epsilon > 0$ we consider the ($\epsilon$-discrete)
metric $\min(d,\epsilon)$ on $X$ and ($\epsilon$-bounded) metric
$\max(d,\epsilon)$~\cite{Brod-Dydak-Higes-Mitra}. Therefore we can
define and work with all-scales notions and then obtain the large
scale (or small scale) results as corollaries after applying the
macroscopic (or microscopic) functor.

In this paper we consider five categories of separable metric
spaces: Lipschitz, Uniform, the corresponding Metrically Proper
subcategories (see the definitions at the end of Introduction),
and the Coarse category defined by Roe~\cite{Roe lectures}.

The concept of dimension appropriate for the Lipschitz category is
the Assouad-Nagata dimension~\cite{Lang-Sch Nagata dim}. For
discrete metric spaces the notion of Assouad-Nagata dimension is
equivalent to the notion of asymptotic dimension of linear type
considered by Gromov~\cite{Gro asym invar} and Roe~\cite{Roe
lectures} (Dranishnikov and Zarichnyi call it "asymptotic
dimension with Higson property"~\cite{Dran-Zar}). For bounded
metric spaces the notion of Assouad-Nagata dimension is equivalent
to the notion of capacity dimension introduced recently by
Buyalo~\cite{Buyalo1, Buyalo Lebedeva}.

In Section~\ref{section uniform} we introduce the concept of
dimension appropriate for the uniform category. For discrete metric
spaces the notion of uniform dimension is equivalent to the notion
of asymptotic dimension introduced by Gromov. For a bounded metric
space $X$ the uniform dimension $\uudim(X)$ coincides with the
large dimension $\Delta d X$ from the book~\cite{Isbell book}.

Ultrametric spaces play the central role in this paper. We show
that in (Proper) Lipschitz and (Proper) Uniform categories a
metric space $(X,d)$ has dimension 0 if and only if there is an
ultrametric $\rho$ on $X$ such that the identity map $(X,d)\to
(X,\rho)$ is an equivalence (for separable metric spaces and
continuous maps this result was proved by de~Groot~\cite{Groot}
and Nagata~\cite{Nagata}; for metric spaces and Lipschitz maps it
is proved in~\cite[Chapter 15]{DavSem}; for discrete spaces and
coarse maps this result belongs to M.~Zarichnyi~\cite{MZar
commun}). We also exhibit an ultrametric space which is universal
(in all categories) for all 0-dimensional spaces. Notice that
there is an ultrametric space containing isometric copy of any
ultrametric space~\cite{Lemin Lemin universal, Lemin
ultrametrization, Bogatyi 2}.

In (Proper) Lipschitz and (Proper) Uniform categories we
characterize 0-dimensional spaces by means of retractions to
subspaces. In the Lipschitz category we prove that the following
conditions are equivalent:
\begin{enumerate}
\item $X$ has dimension 0; \item the unit 0-sphere $S^0$ is an
absolute extensor for $X$; \item every metric space is an absolute
extensor for $X$.
\end{enumerate}
We failed to find the analogous characterization in the Uniform
category.

In Sections~\ref{section groups} and~\ref{section Examples} we
consider discrete metric spaces in the Coarse category. It is easy
to see that a finitely generated group of asymptotic dimension 0
is finite and therefore all such groups are coarsely equivalent.
To define asymptotic dimension for an infinitely generated
countable group one should consider a left invariant proper metric
on it \cite{Sha}. We describe a natural way to introduce such a metric and
prove that any group of asymptotic dimension 0 is coarsely
equivalent to an abelian group. It is known that a countable group
has asymptotic dimension 0 if and only if it is locally
finite~\cite{Smith} but we are not aware of any characterization
of locally finite countable groups up to coarse equivalence. In
Section~\ref{section Examples} we construct uncountably many
examples of coarsely inequivalent metric spaces of asymptotic
dimension 0. The idea of the construction does not work for
groups.

\begin{Def}
A map $f\colon X\to Y$ of metric spaces is called {\it Lipschitz}
if there is a constant $\lambda>0$ such that the inequality
$d_Y(f(x),f(x'))\le\lambda\cdot d_X(x,x')$ holds for all points
$x,x'\in X$. $f$ is called {\it $\lambda$-Lipschitz} if we need to
specify the constant $\lambda$. $f$ is called {\it
$\lambda$-bi-Lipschitz} if both $f$ and $f^{-1}$ are
$\lambda$-Lipschitz.
\end{Def}

For any Lipschitz map $f$ we denote
$$ Lip(f)=\inf\{\lambda\mid \text{$f$ is $\lambda$-Lipschitz}\} $$
Notice that a Lipschitz map $f$ is $Lip(f)$-Lipschitz.

\begin{Def}
A metric space $X$ is called a {\it Lipschitz extensor for a metric
space $Y$} if there exists a constant $m>0$ such that for any
closed subspace $A\subset Y$ any Lipschitz map $f\colon A\to X$
extends to a Lipschitz map $F\colon Y\to X$ with $Lip(F)=m\times
Lip(f)$. We call the space $X$ an {\it $m\times$-Lipschitz extensor
for $Y$} if we need to specify the constant $m$.
\end{Def}

A map $f\colon X\to Y$ is called {\it metrically proper} if for
any bounded subset $A$ of the space $Y$ its preimage $f^{-1}(A)$
is bounded.

\begin{Def}
The {\it Lipschitz category} consists of separable metric spaces with
morphisms being Lipschitz maps. Its subcategory of unbounded
spaces and metrically proper maps is called the {\it Proper Lipschitz
category}.
\end{Def}

We call a map $f\colon X\to Y$ {\it uniform} if there is a
function $\delta_f\colon \RR_+\to \RR_+$ with $\lim_{t\to
0}\delta_f(t)=0$ such that $d_Y(f(x),f(x'))\le\delta_f(d_X(x,x'))$
for all points $x,x'\in X$. To specify the function $\delta_f$ we
sometimes say that the map $f$ is $\delta_f$-uniform. A map $f$ is
called {\it bi-uniform} if both $f$ and $f^{-1}$ are uniform.

\begin{Def}
The {\it Uniform category} consists of separable metric spaces with
morphisms being uniform maps. Its subcategory of unbounded spaces
and metrically proper maps is called the {\it Proper Uniform
category}.
\end{Def}

We call a metric space $X$ {\it discrete} if there is $\epsilon >
0$ such that $X$ is $\epsilon$-discrete.

We call a map $f\colon X\to Y$ {\it large scale uniform} if there
is a function $\delta_f\colon \RR_+\to \RR_+$ such that
$d_Y(f(x),f(x'))\le\delta_f(d_X(x,x'))$ for all points $x,x'\in
X$. A map is called {\it coarse} if it is large scale uniform and
metrically proper. Metric spaces $X$ and $Y$ are {\it coarsely
equivalent} if there exist a constant $C>0$ and two coarse maps
$f\colon X\to Y$ and $g\colon Y\to X$ such that the maps $g\circ
f$ and $f\circ g$ are $C$-close to the identity.

\section{Ultrametric spaces}\label{section Ultrametric spaces}

\begin{Def}
A metric space $(X,d)$ is called {\it ultrametric} if for all
$x,y,z\in X$ we have $d(x,z) \le \max\{d(x,y),d(y,z)\}$.
\end{Def}

An ultrametric space $X$ can be characterized by the following
very useful property:

{\bf Ultametric property of a triangle}. If a triangle in a space
$X$ has sides (distances between vertices) $a\le b\le c$, then
$b=c$.

The following properties of ultrametric space are easy to check. A
ball of radius $D$ in an ultrametric space has diameter $D$. Two
balls of radius $D$ in an ultrametric space are either
$D$-disjoint or identical.

\begin{Prop}\label{ultrametric characterization}
Let $(X,d)$ be a metric space. The metric $d$ is an ultrametric if
and only if $f(d)$ is a metric for every nondecreasing function
$f\colon \RR_+\to\RR_+$.
\end{Prop}

\begin{pf}
If $d$ is ultrametric and $a\le b = c$ are sides of a triangle in
$(X,d)$ then $f(a)\le f(b) = f(c)$ are sides of the corresponding
triangle in $(X,f(d))$ and therefore $f(d)$ is an ultrametric.

If $d$ is not an ultrametric then there is a triangle in $(X,d)$
with sides $a\le b < c$. Consider the function
$$ f(t)\begin{cases} t\text{ if } t\le b \\
\frac{2b}{c-b}t+\frac{bc-3b^2}{c-b}\text{ if } t\ge b
\end{cases}
$$
The sides of the corresponding triangle in $(X,f(d))$ are
$f(a)\le f(b)=b < 3b= f(c)$ which contradicts the triangle
inequality.
\end{pf}

\begin{Def}
A metric is said to be {\it $3^n$-valued} if the only values assumed by the
metric are $3^n$, $n\in\ZZ$.
\end{Def}

The triangle inequality for a metric $d$ implies the following:

\begin{Lem}\label{3^n-valued metric is ultrametric}
Any $3^n$-valued metric is an ultrametric.
\end{Lem}

\begin{Lem}\label{bi-Lipschitz change to 3^n-valued ultrametric}
Any ultrametric space is 3-bi-Lipschitz equivalent to a
$3^n$-valued ultrametric space.
\end{Lem}

\begin{pf}
Given an ultrametric space $(X,d)$ we define a new metric $\rho$
on $X$ as follows:
$$\rho(x,y)=3^n \text{\qquad if \qquad} 3^{n-1}<d(x,y)\le 3^n. $$
Clearly, the identity map $\id\colon (X,d)\to (X,\rho)$ is
expanding and 3-Lipschitz.
\end{pf}

Let us describe an ultrametric space $(\LL_\omega,\mu)$ which is
universal for all separable ultrametric spaces with $3^n$-valued
metrics. This space appeared naturally in different areas of
mathematics (see for example~\cite{Lemin ultrametrization} and
references therein). Let us fix a countable set $S$ with a
distinguished element $s_0\in S$. The set $\LL_\omega$ is a subset
of the set of infinite sequences $\bar{x}=\{x_n\}_{n\in \ZZ}$ with
all elements $x_n$ from the set $S$. A sequence
$\bar{x}$ belongs to $\LL_\omega$ if there exists an index $k\in
\ZZ$ such that $x_n=s_0$ for all $n<k$. The metric $\mu$ is
defined as $\mu(\bar{x},\bar{y})=3^{-m}$ where $m\in\ZZ$ is the
minimal index such that $x_m\ne y_m$. Clearly, the space
$\LL_\omega$ is a complete separable ultrametric (by
Lemma~\ref{3^n-valued metric is ultrametric}) space.

To prove that any separable ultrametric space with $3^n$-valued
metric embeds isometrically into $(\LL_\omega,\mu)$ we follow the
idea of P.S.~Urysohn~\cite{Urysohn} and show that the space
$\LL_\omega$ is {\it finitely injective}:

\begin{Lem}\label{L omega is finitely injective}
Let $(X,d)$ be a finite metric space with $3^n$-valued metric $d$.
For any subspace $A\subset X$, any isometric map $f\colon A\to
\LL_\omega$ admits an isometric extension $\tilde{f}\colon X\to
\LL_\omega$.
\end{Lem}

\begin{pf}
It is sufficient to prove Lemma in case $X\setminus A$ consists of one
point $x$. In such case we have to find a point $\bar{z}\in
\LL_\omega$ such that $\mu(\bar{z},f(a))=d(x,a)$ for every point
$a\in A$. Let $A_x=\{a\in A\mid d(x,a)=d(x,A)\}$ be the set of all
points in $A$ closest to $x$ and let $d(x,A)=3^{-n}$. Fix a point
$b\in A_x$ and define $\bar{z}=\{z_n\}_{n\in \ZZ}$ as follows:
$z_m=f(b)_m$ if $m<n$; $z_m=s_0$ if $m>n$; $z_n$ is any element of
the set $S$ other than $f(c)_n$ for any point $c\in A_x$.

Clearly, $\mu(\bar{z},f(c))=3^{-n}=d(x,c)$ for any point $c\in
A_x$. For any point $a\in A\setminus A_x$ we have
$d(a,x)=d(a,b)=3^{-m}>3^{-n}$ which means that $f(a)_m\ne
f(b)_m=z_m$ and therefore $\mu(\bar{z},f(a))=3^{-m}=d(x,a)$.
\end{pf}

\begin{Thm}\label{L omega is isometry-universal}
Any separable metric space $(X,d)$ equipped with $3^n$-valued metric $d$
embeds isometrically into the space $(\LL_\omega,\mu)$.
\end{Thm}

\begin{pf}
Since $X$ is separable, it is sufficient to embed isometrically a
countable dense subspace $A$ of $X$. One can embed such a subspace
by induction using Lemma~\ref{L omega is finitely injective}.
\end{pf}

\begin{Cor}\label{L omega is bi-Lipschitz universal for ultrametric spaces}
Any separable ultrametric space admits 3-bi-Lipschitz embedding
into the space $(\LL_\omega,\mu)$.
\end{Cor}

\begin{pf}
Combine Lemma~\ref{bi-Lipschitz change to 3^n-valued ultrametric}
and Theorem~\ref{L omega is isometry-universal}.
\end{pf}

\begin{Thm}\label{UnboundedIsRetract}
Every closed subset $A$ of an ultrametric space $X$ is a
$\lambda$-Lipschitz retract of $X$ for any $\lambda>1$. If the
subset $A$ is unbounded, the retraction can be chosen to be
metrically proper.
\end{Thm}

\begin{pf}
Suppose that $X$ is an ultrametric space and $A\subset X$ is a
closed subspace. If $\lambda>1$ is given, choose a number
$\delta>1$ such that $\delta^2<\lambda$.

Let us fix a base point $x_0\in X$. Take an arbitrary well-order
$<_{k}$ on each non empty Annulus $A_k = \{x | k \le d(x, x_0) <
k+1\}$ of $X$ for every $k \in \NN \cup \{0\}$. Now we say $z
\prec z'$ for any two points $z, z' \in X$ if $z \in A_k$, $z' \in
A_{k'}$ and $k > k'$ or if $z, z' \in A_k$ and $z <_{k} z'$.
Notice that $\prec$ is an order in $X$ such that for every non
empty bounded subset $C$ of $X$ the restricted order $\prec|_{C}$
is a well-order.

We define a retraction $r\colon X\to A$ as follows. For a point
$x\in X$ we look at the nonempty bounded set
$$A_x=\{a\in A\mid d(x,a)\le\delta\cdot\dist(x,A)\}$$
and put $r(x)$ to be the minimal point in the set $A_x$ with
respect to the order $\prec$.

Let us show that the retraction $r$ is $\lambda$-Lipschitz. Assume
that for some points $x,y\in X$ we have $d(r(x),r(y))>\lambda\cdot
d(x,y)$. Without loss of generality we may assume that $r(x)\prec
r(y)$.

If $d(y,r(x))\le d(y,r(y))$, then $r(x)\in A_y$ and $r(x)\prec
r(y)$ contradicts the choice of $r(y)$ to be the minimal point in
the set $A_y$.

In case $d(y,r(x))>d(y,r(y))$ we denote by $D$ the distance
between $r(x)$ and $r(y)$ and notice that
$d(y,r(x))=d(r(x),r(y))=D$ in the isosceles triangle $\{y, r(x),
r(y)\}$. Since $D>d(x,y)$, we have $d(x,r(x))=d(y,r(x))=D$ in the
isosceles triangle $\{x, y, r(x)\}$.

$$d(x,r(y))\ge \dist(x,A)\ge \frac{1}{\delta}\cdot
d(x,r(x))=\frac{D}{\delta}>\frac{D}{\lambda}>d(x,y) $$

Therefore $d(x,r(y))=d(y,r(y))$ in the isosceles triangle $\{x, y,
r(y)\}$. The point $r(x)$ does not belong to $A_y$ since
$r(x)\prec r(y)$, thus $d(y,r(x))=D>\delta\cdot\dist(y,A)$. Then
there exists a point $z\in A$ with $d(y,z)<\frac{D}{\delta}$.

$$d(y,z)\ge\dist(y,A)\ge\frac{d(y,r(y))}{\delta}=
\frac{d(x,r(y))}{\delta}\ge\frac{D}{\delta^2}>\frac{D}{\lambda}>d(x,y)
$$

Therefore $d(x,z)=d(y,z)$ in the isosceles triangle $\{x, y, z\}$.
Since $d(x,z)<d(x,r(x))$, we have $z\in A_x$, but
$d(x,z)<\frac{D}{\delta}=\frac{d(x,r(x))}{\delta}$ contradicts the
definition of $A_x$ (two points $a,a'\in A_x$ cannot satisfy
$d(x,a)<\frac{d(x,a')}{\delta}$).

If the subset $A$ is unbounded, we prove that the retraction $r$
is metrically proper. Let $B$ be any bounded subset of $A$. Choose
a point $a \in A$ which is in an annulus greater than any annulus
that has non-empty intersection with $B$ (therefore, $a\prec B$).
Given any point $x \in r^{-1}(B)$ we have $a\not\in A_x$,
therefore $d(x, r(x)) \le \delta \cdot d(x, A) < d(x, a)$. The
ultrametric property of the triangle $\{x, a, r(x)\}$ implies
$d(r(x), a) = d(x, a)$ therefore:
$$ d(x , B) \le d(x, r(x)) < d(r(x), a) \le diam(B)+ d(a, B)$$
\end{pf}

\begin{ex}
Let $X=\{x_n\}_{n=1}^\infty$ be a sequence of points. Define
$d(x_1,x_n)=1+\frac{1}{n}$ and
$d(x_m,x_n)=\max\{1+\frac{1}{m},1+\frac{1}{n}\}$ for any $m,n>1$.
Then $d$ is an ultrametric on $X$ and there is no 1-Lipschitz
retraction of $X$ onto $A=\{x_n\}_{n=2}^\infty$.
\end{ex}

\section{Assouad-Nagata dimension}\label{section Nagata}

\begin{Def}
Let $X$ be a metric space, $A$ be a subspace of $X$, and $S$ be a
positive number.

$A$ is {\it $S$-bounded} if for any points $x,x'\in A$ we have
$d_X(x,x')\le S$.

An {\it $S$-chain} in $A$ is a sequence of points $x_1,\ldots,x_k$
in $A$ such that for every $i<k$ the set $\{x_i,x_{i+1}\}$ is
$S$-bounded.

$A$ is {\it $S$-connected} if for any points $x,x'\in A$ can be
connected in $A$ by an $S$-chain.
\end{Def}

Notice that any subset $A$ of $X$ is a union of its {\it
$S$-components} (the maximal $S$-connected subsets of $A$). If $B$
and $B'$ are two $S$-components of the set $A$ then $B$ and $B'$
are $S$-disjoint. Intuitively, a metric space $X$ has dimension
$0$ at scale $S>0$ if all $S$-components of $X$ are uniformly
bounded.

\begin{Def}
A metric space $X$ has Assouad-Nagata dimension zero (notation
$dim_{AN}(X)\le 0$) if there exists a constant $m\ge 1$, such that
for any $S>0$ all $S$-components of $X$ are $mS$-bounded.
\end{Def}

It is easy to see that bi-Lipschitz maps preserve Assouad-Nagata
dimension.

Ultrametric spaces are the best examples of metric spaces of
Assouad-Nagata dimension zero. Indeed, for any positive number $D$
any $D$-component of an ultrametric space is a $D$-ball and
therefore is $D$-bounded. Let us characterize spaces of
Assouad-Nagata dimension 0 using ultrametrics.

The following theorem is proved in~\cite[Proposition
15.7]{DavSem}. We provide a proof for completeness.

\begin{Thm}\label{Nagata dim vs. ultrametrics}
If a metric space $(X,d)$ has Assouad-Nagata dimension
$\dim_{AN}(X)\le 0$, then there is an ultrametric $\rho$ on $X$
such that the identity map $\id\colon (X,d)\to (X,\rho)$ is
bi-Lipschitz.
\end{Thm}

\begin{pf}
Suppose that for a number $m>1$, all $S$-components of $X$ are
$mS$-bounded. Consider two points $x,z\in X$ and put
$$S=\frac{d(x,z)}{2m}.$$
Then the points $x$ and $z$ belong to different $S$-components of
$X$. Thus for any chain $x=x_0,x_1,\ldots,x_{k-1}, x_k=z$ we have
$$ d(x,z)\le 2m\cdot\max_{0\le i<k}\{d(x_i,x_{i+1})\}. $$

Now define $\rho(x,z)$ to be the infimum of $\max_{0\le
i<k}\{d(x_i,x_{i+1})\}$ over all finite chains
$x_0,x_1,\ldots,x_{k-1}, x_k$ with $x=x_0$ and $x_k=z$. Clearly
$$\frac{1}{2m}\cdot d(x,z)\le \rho(x,z)\le d(x,z). $$

To see that $\rho$ is an ultrametric, take three points $x,y,z$ in
$X$ and let $s$ be the infimum of all positive numbers $S$ such
that all three points belong to one $S$-component of $X$. If all
three points belong to one $s$-component or all three belong to
different $s$-components, then $\rho(x,y)=\rho(x,z)=\rho(y,z)=s$.
If the points $x$ and $y$ belong to one $s$-component which does
not contain $z$, then $\rho(x,y)\le s=\rho(x,z)=\rho(y,z)$.
\end{pf}

\begin{Thm}\label{Universal space for Nagata dim 0}
Any separable metric space of Assouad-Nagata dimension 0 admits a
bi-Lipschitz embedding into the space $(\LL_\omega,\mu)$.
\end{Thm}

\begin{pf}
Apply Theorem~\ref{Nagata dim vs. ultrametrics} and Theorem~\ref{L
omega is bi-Lipschitz universal for ultrametric spaces}.
\end{pf}

\begin{Thm}\label{Nagata dim characterization}
In the Lipschitz category the following conditions are equivalent:
\begin{enumerate}
\item $\dim_{AN}(X)\le 0$; \item there exists a number $\lambda$
such that every closed subset of $X$ is a $\lambda$-Lipschitz
retract of $X$; \item there exists a number $\lambda$ such that
every metric space is a $\lambda\times$-Lipschitz extensor for
$X$; \item the unit 0-sphere $S^0$ is an extensor for $X$.
\end{enumerate}
Conditions $(1)$, $(2)$, and $(3)$ are equivalent in the Proper
Lipschitz category.
\end{Thm}

\begin{pf}
$(1)\implies (2)$ in both Lipschitz and Proper Lipschitz
categories. Theorem~\ref{Nagata dim vs. ultrametrics} allows us to
find an ultrametric $\rho$ on $X$ which is bi-Lipschitz equivalent
to $d$. Application of Theorem~\ref{UnboundedIsRetract} completes
the proof.

$(2)\implies (3)$ in both Lipschitz and Proper Lipschitz
categories. Given a closed subspace $A\subset X$ and a Lipschitz
map $f\colon A\to Y$ to some metric space $Y$ we fix a
$\lambda$-Lipschitz retraction $r\colon X\to A$. Then the
composition $f\circ r\colon X\to K$ has the Lipschitz constant
bounded by $\lambda\cdot\Lip(f)$.

$(3)\implies (4)$ Obvious.

$(4)\implies (1)$ Let $m\ge 1$ be a number such that any
$\lambda$-Lipschitz map from any closed subspace $A\subset X$ to
$S^0$ can be extended to $m\lambda$-Lipschitz map of $X$. If an
$S$-component of $X$ is not $mS$-bounded, there are points $z_0$
and $z_1$ with $d(z_0,z_1)>mS$ and an $S$-chain of points
$z_0=x_0,x_1,\ldots,x_k=z_1$. Notice that the map $f\colon
\{z_0\}\cup\{z_1\}\to S^0$ defined as $f(z_0)=0$ and $f(z_1)=1$ is
$\frac{1}{d(z_0,z_1)}$-Lipschitz but any extension of this map to
the chain is at least $\frac{1}{S}$-Lipschitz and cannot be
$\frac{m}{d(z_0,z_1)}$-Lipschitz (since
$\frac{1}{S}>\frac{m}{d(z_0,z_1)}$).

$(3)\implies (1)$ in the Proper Lipschitz category. If an
$S$-component of $X$ is not $\lambda S$-bounded, there are points
$z_0$ and $z_1$ with $d(z_0,z_1)>\lambda S$ and an $S$-chain of
points $z_0=x_0,x_1,\ldots,x_k=z_1$. Let $A$ be any unbounded
$\lambda S$-discrete subspace of $X$ containing the points $z_0$
and $z_1$. Notice that the identity map $\id_A$ is $1$-Lipschitz
but any extension of this map to the chain is not $\lambda
S$-Lipschitz.
\end{pf}

\begin{Problem}
Is there an analog of condition $(4)$ from Theorem~\ref{Nagata dim
characterization} in the Proper Lipschitz category?
\end{Problem}

\section{Uniform dimension}\label{section uniform}

\begin{Def}\label{uniform dimension}
A metric space $X$ has {\it uniform} dimension zero (notation
$dim_u(X)\le 0$) if there exists a continuous increasing function
$\uD\colon \RR_+\to\RR_+$ with $\uD(0)=0$ and $\lim_{t\to
\infty}\uD(t)=\infty$, such that for any positive number $S$ every
$S$-component of $X$ is $\uD(S)$-bounded.

To specify the function $\uD$ we sometimes say that the space $X$
has {\it uniform dimension zero of type $\uD$}.
\end{Def}

If the function $\uD$ does not exceed some linear function
$\uD(t)\le k\cdot t$ for all $t\ge 0$, then the space $X$ has
Assouad-Nagata dimension 0. We want the dimension control function
to be increasing and continuous to guarantee the existence of the
inverse function $\uD^{(-1)}$.

It is easy to check that the uniform dimension is preserved under
the bi-uniform maps:

\begin{Lem}\label{uniform dim under uniform map}
Let $f\colon X\to Y$ be a bi-uniform map. Then
$\dim_u(X)=\dim_u(f(X))$.
\end{Lem}

\begin{Thm}\label{uniform dim vs. ultrametrics}
If a metric space $(X,d)$ has uniform dimension $\dim_u(X)\le 0$,
then there is an ultrametric $\rho$ on $X$ such that the identity
map $\id\colon (X,d)\to (X,\rho)$ is bi-uniform.
\end{Thm}

\begin{pf}
Suppose that the space $X$ has uniform dimension zero of type
$\uD$. Consider two points $x,z\in X$ and put
$$S=\frac{1}{2}\uD^{-1}(d(x,z)).$$
Then the points $x$ and $z$ belong to different $S$-components of
$X$. Thus for any chain $x=x_0,x_1,\ldots,x_{k-1}, x_k=z$ we have
$$ \uD^{-1}(d(x,z))\le 2\cdot\max_{0\le i<k}\{d(x_i,x_{i+1})\}. $$

Now define $\rho(x,z)$ to be the infimum of $\max_{0\le
i<k}\{d(x_i,x_{i+1})\}$ over all finite chains
$x_0,x_1,\ldots,x_{k-1}, x_k$ with $x=x_0$ and $x_k=z$. Clearly
$$\frac{1}{2}\cdot \uD^{-1}(d(x,z))\le \rho(x,z)\le d(x,z). $$

To see that $\rho$ is an ultrametric, take three points $x,y,z$ in
$X$ and let $s$ be the infimum of all positive numbers $S$ such
that all three points belong to one $S$-component of $X$. If all
three points belong to one $s$-component or all three belong to
different $s$-components, then $\rho(x,y)=\rho(x,z)=\rho(y,z)=s$.
If the points $x$ and $y$ belong to one $s$-component which does
not contain $z$, then $\rho(x,y)\le s=\rho(x,z)=\rho(y,z)$.
\end{pf}

\begin{Cor}\label{L omega is universal for uniform dimension}
A separable metric space $X$ has uniform dimension zero if and
only if it admits a bi-uniform embedding into $\LL_\omega$.
\end{Cor}

\begin{pf}
If $\dim_u(X)\le 0$ we can change the metric on $X$ bi-uniformly
to get an ultrametric space and then embed it in a bi-Lipschitz
way into $\LL_\omega$ using Theorem~\ref{L omega is bi-Lipschitz
universal for ultrametric spaces}.

If $X$ embeds bi-uniformly into $\LL_\omega$, its image has
uniform dimension zero as a subspace of $\LL_\omega$. Then $X$ has
uniform dimension zero by Lemma~\ref{uniform dim under uniform
map}.
\end{pf}

\begin{Thm}\label{Uniform dim characterization}
In both Uniform and Proper Uniform categories the following
conditions are equivalent:
\begin{enumerate}
\item $\uudim X\le 0$; \item there exists a continuous increasing
function $\mu\colon \RR_+\to\RR_+$ with $\mu(0)=0$ and $\lim_{t\to
\infty}\mu(t)=\infty$, such that every closed subspace of $X$ is
$\mu$-uniform retract of $X$.
\end{enumerate}
\end{Thm}

\begin{pf}
$(1)\implies (2)$ Theorem~\ref{uniform dim vs. ultrametrics}
allows us to find an ultrametric $\rho$ on $X$ which is
bi-uniformly equivalent to $d$. Application of
Theorem~\ref{UnboundedIsRetract} completes the proof.

$(2)\implies (1)$ If an $S$-component of $X$ is not
$\mu(S)$-bounded, there are points $z_0$ and $z_1$ with
$d(z_0,z_1)>\mu(S)$ and an $S$-chain of points
$z_0=x_0,x_1,\ldots,x_k=z_1$.

In the Uniform category let $A=\{z_0\}\cup\{z_1\}$. In the Proper
Uniform category we consider any unbounded closed subspace $A$ of
$X$ containing the points $z_0$ and $z_1$ and such that the
distance from $\{z_0\}\cup\{z_1\}$ to the rest of $A$ is greater
than $d(z_0,z_1)$.

Notice that any retraction of $X$ onto $A$ restricted to the chain
takes some $S$-closed points to two points of distance greater
than $d(z_0,z_1)>\mu(S)$. Thus such a retraction cannot be
$\mu$-uniform.
\end{pf}

\begin{Problem}
Are there analogs of conditions $(3)$ and $(4)$ from
Theorem~\ref{Nagata dim characterization} in the Proper Uniform
category?
\end{Problem}

\section{Locally finite countable groups}\label{section groups}

It is proved in~\cite{Smith} that a countable group $G$ (equipped
with any proper metric) has asymptotic dimension zero if and only
if $G$ is locally finite (i.e. every finitely generated subgroup
of $G$ is finite). The purpose of this section is to show that
such a group is bi-uniformly equivalent to a locally finite
abelian group. Also we classify locally finite countable groups up
to bi-uniform equivalence. The problem of classification of
locally finite countable groups up to coarse equivalence remains
open. Notice that for discrete metric spaces the notions of
bi-uniform equivalence and bijective coarse equivalence coincide.

A left invariant metric $d$ on a countable group $G$ is {\it
proper} if and only if every bounded subset of $(G,d)$ is finite.
Thus a left invariant proper metric $d$ on $G$ is bounded from
below and therefore the asymptotic dimension of $(G,d)$ is equal
to its uniform dimension. There is only one way (up to bi-uniform
equivalence) to introduce a proper left-invariant metric on
$G$~\cite[Proposition 1]{Smith}. Thus the asymptotic dimension of
a countable group does not depend on the choice of a proper
left-invariant metric.

Let $G$ be a locally finite countable group. Let us describe a
particularly simple way to define a proper left-invariant metric
on $G$. Consider a filtration ${\mathcal{L}}$ of $G$ by finite
subgroups ${\mathcal{L}} =\{1\subset G_1\subset G_2 \subset G_3
\dots\}$ and define the metric $d_{\mathcal{L}}$ associated to
this filtration as:
$$d_{\mathcal{L}}(x,y)= \min\{i \mid x^{-1}y\in G_i \}.$$

Clearly, $d_{\mathcal{L}}$ is an ultrametric (therefore, the
asymptotic dimension of $(G,d_{\mathcal{L}})$ is zero).

\begin{Lem}\label{IsometricLocFinGroups}
Suppose two groups $G$ and $H$ have filtrations by finite
subgroups: ${\mathcal{L}} =\{1\subset G_1\subset G_2 \subset G_3
\dots\}$ of $G$ and ${\mathcal{K}} =\{1\subset H_1\subset H_2
\subset H_3 \dots\}$ of $H$. If the index $[G_i:G_{i-1}]$ is less
than or equal to the index $[H_i:H_{i-1}]$ for all $i$, then
$(G,d_{\mathcal{L}})$ admits an isometric embedding into
$(H,d_{\mathcal{H}})$. Moreover, if
$[G_i:G_{i-1}]=[H_i:H_{i-1}]$ for all $i$
(equivalently, the cardinality of $G_i$ equals cardinality of $H_i$ for all $i$), then the groups
$(G,d_{\mathcal{L}})$ and $(H,d_{\mathcal{H}})$ are isometric.
\end{Lem}

\begin{pf}
Put $a_i=[G_i:G_{i-1}]$ and $b_i=[H_i:H_{i-1}]$. Fix an injection
$f_1\colon G_1\to H_1$ and assume injections $f_k\colon G_k\to
H_k$ are known for $k\leq n$ such that the following two
properties hold:
 \begin{enumerate}
\item  $f_i(x)=f_j(x)$ for $i < j$ and $x\in G_i$, \item the
injection $f_k\colon G_k\to H_k$ is isometric.
\end{enumerate}

Pick an injection of the set of cosets $\{x\cdot G_n\}$ of $G_n$
in $G_{n+1}$ into the set of cosets $\{y\cdot H_n\}$ of $H_n$ in
$H_{n+1}$. That amounts to picking representatives $1$,
$x_1$,\ldots, $x_m$ ($m=a_{n+1}-1$) of cosets of $G_n$ in
$G_{n+1}$ and picking representatives $1$, $y_1$,\ldots, $y_l$
($l=b_{n+1}-1$) of cosets of $H_n$ in $H_{n+1}$. Make sure the
injection takes $\{1\cdot G_n\}$ to $\{1\cdot H_n\}$. Now we
extend $f_n$ to $f_{n+1}\colon G_{n+1}\to H_{n+1}$ as follows: if
$x\in G_{n+1}\setminus G_n$, we represent $x$ as $x_k\cdot x'$ for
some unique $k\leq m$ and we define $f_{n+1}(x)$ as $y_k\cdot
f_n(x')$.

If $x$ and $z$ belong to different cosets of $G_n$ in $G_{n+1}$,
then $f_{n+1}(x)$ and $f_{n+1}(z)$ belong to different cosets of
$H_n$ in $H_{n+1}$ and
$d_{\mathcal{L}}(x,z)=n+1=d_{\mathcal{H}}(f_{n+1}(x),f_{n+1}(z))$.
If $x$ and $z$ belong to the same coset $x_k\cdot G_n$ of $G_n$ in
$G_{n+1}$, then $x=x_k\cdot x'$, $z=x_k\cdot z'$. Since
$f_{n+1}(x)=y_k\cdot f_n(x')$, $f_{n+1}(z)=y_k\cdot f_n(z')$, and
the map $f_n$ is isometry, then
$$d_{\mathcal{L}}(x,z)=d_{\mathcal{L}}(x',z')d_{\mathcal{H}}(f_{n}(x'),f_{n}(z'))d_{\mathcal{H}}(f_{n+1}(x),f_{n+1}(z)).$$

By pasting all $f_n$ we get an isometric injection $f\colon G\to
H$. Notice that in case $[G_i:G_{i-1}]=[H_i:H_{i-1}]$ for all $i$,
the map $f$ is bijective and establishes an isometry between
$(G,d_{\mathcal{L}})$ and $(H,d_{\mathcal{H}})$.
\end{pf}

\begin{Lem}\label{CharOfIsometricLocFinGroups}
Given two locally finite groups $G$ and $H$
the following conditions are equivalent:
\begin{enumerate}
\item There are left-invariant proper metrics $d_G$ on $G$ and $d_H$
on $H$ such that $(G,d_G)$ is isometric to $(H,d_H)$.
\item There are filtrations by finite
subgroups: ${\mathcal{L}} =\{1\subset G_1\subset G_2 \subset G_3
\dots\}$ of $G$ and ${\mathcal{K}} =\{1\subset H_1\subset H_2
\subset H_3 \dots\}$ of $H$ such that the cardinality of $G_i$ equals cardinality of $H_i$ for all $i$.
\end{enumerate}
\end{Lem}

\begin{pf}
In view of \ref{IsometricLocFinGroups}, it suffices to show (1)$\implies$(2).
Obviously, we may pick an isometry $f\colon G\to H$ such that $f(1_G)=1_H$
(replace any $f$ by $f(1_G)^{-1}\cdot f$).
Notice $f$ establishes bijectivity between $m$-component of $G$ containing
$1_G$ and the $m$-component of $H$ containing $1_H$. Also, those components are
subgroups of $G$ and $H$. Thus, define $G_1$ as $1$-component of $G$ containing
$1_G$ and, inductively, $G_{i+1}$ as $(\diam(G_i)+i)$-component
of $G$ containing $1_G$.
\end{pf}

{\bf Main example}. If $G$ is a direct sum of cyclic groups
$\bigoplus\limits_{i=1}^\infty {\mathbb{Z}}_{a_i}$ we consider the metric
on $G$ associated to the filtration $${\mathcal{L}}
=\{1\subset{\mathbb{Z}}_{a_1}\subset
{\mathbb{Z}}_{a_1}\oplus{\mathbb{Z}}_{a_2}\subset
{\mathbb{Z}}_{a_1}\oplus{\mathbb{Z}}_{a_2}\oplus{\mathbb{Z}}_{a_3}\subset\dots\}$$

If we write elements of the group $\bigoplus\limits_{i=1}^\infty
{\mathbb{Z}}_{a_i}$ as $p=p_1p_2\dots p_n$ where
$p_j\in{\mathbb{Z}}_{a_j}$ and denote $|p|=n$ then the ultrametric
$d_{\mathcal{L}}$ can be defined explicitly as
$$
d_{\mathcal{L}}(p,q)=\begin{cases} \max\{|p|,|q|\}\text{ if }
|p|\ne |q|\\
\max\{i\mid p_i\ne q_i\}\text{ if } |p|=|q|
\end{cases}
$$

\begin{Thm}\label{MetricOnLocallyFiniteGrps}
A locally finite countable group $G$ with a proper left invariant
metric $d$ is bi-uniformly equivalent to a direct sum of cyclic
groups.
\end{Thm}

\begin{pf}
Fix a filtration ${\mathcal{L}}$ of $G$ by finite subgroups
${\mathcal{L}} =\{1\subset G_1\subset G_2 \subset G_3 \dots\}$.
Then $(G,d)$ is bi-uniformly equivalent to
$(G,d_{\mathcal{L}})$~\cite[Proposition 1]{Smith}. By
\ref{IsometricLocFinGroups}, $(G,d_{\mathcal{L}})$ is isometric to
$\bigoplus\limits_{i=1}^\infty {\mathbb{Z}}_{a_i}$ where
$a_i=[G_i:G_{i-1}]$.
\end{pf}

\begin{Def}
Let $G$ be a countable locally finite group and $p$ be a prime
number. We define a {\it $p$-Sylow number of $G$} (finite or
infinite) as follows:
$$ |\text{$p$\,-Syl}|(G)= \sup\{p^n\mid p^n \text{ divides } |F|, F \text{ a finite
subgroup of }G, n\in\ZZ\}$$
\end{Def}

Notice that if the $p$-Sylow number of $G$ is finite, it is equal
to the order of a $p$-Sylow subgroup of some finite subgroup of
$G$. For an abelian torsion group $G$ the $p$-Sylow number of $G$
is equal to the order of the $p$-torsion subgroup of $G$.

We are going to use the following theorem of Protasov:

\begin{Thm}[{\cite[Theorem 5]{ProtasovMorphisms}}]\label{Protasov theorem}
Two countable locally finite groups $G$ and $H$ with proper left
invariant metrics are bi-uniformly equivalent if and only if, for
every finite subgroup $F$ of $G$, there exists a finite subgroup
$E$ of $H$ such that $|F|$ is a divisor of $|E|$, and, for every
finite subgroup $E$ of $H$, there exists a finite subgroup $F$ of
$G$ such that $|E|$ is a divisor of $|F|$.
\end{Thm}

\begin{Cor}\label{isomorphic iff bi-uni equiv}
Let $G$ and $H$ be countable direct sums of finite prime cyclic
groups. Let $d_G$ and $d_H$ be proper left invariant metrics on
$G$ and $H$. Then the metric spaces $(G,d_G)$ and $(H,d_H)$ are
bi-uniformly equivalent if and only if the groups $G$ and $H$ are
isomorphic.
\end{Cor}

\begin{Thm}
Let $G$ and $H$ be locally finite countable groups with proper
left invariant metrics $d_G$ and $d_H$. The metric spaces
$(G,d_G)$ and $(H,d_H)$ are bi-uniformly equivalent if and only if
for every prime $p$ we have
$|\text{$p$\,-Syl}|(G)=|\text{$p$\,-Syl}|(H)$.
\end{Thm}

\begin{pf}
Assume the metric spaces $(G,d_G)$ and $(H,d_H)$ are bi-uniformly
equivalent. Our goal is to show that if $|\text{$p$\,-Syl}|(G)\ge
p^n$, then $|\text{$p$\,-Syl}|(H)\ge p^n$. If there is a finite
subgroup $F$ of $G$ such that $p^n$ divides $|F|$, then
by~\ref{Protasov theorem} there is a subgroup $E$ of $H$ such that
$p^n$ divides $|E|$. Thus $|\text{$p$\,-Syl}|(H)\ge p^n$.

Now suppose $|\text{$p$\,-Syl}|(G)=|\text{$p$\,-Syl}|(H)$ for
every prime $p$. By~\ref{Protasov theorem}, it is enough to show
that for every finite subgroup $F$ of $G$, there exists a finite
subgroup $E$ of $H$ such that $|F|$ is a divisor of $|E|$. If
$|F|=p_1^{\alpha_1}\cdot\ldots\cdot p_k^{\alpha_k}$ then
$p_i^{\alpha_i}\le |\text{$p_i$\,-Syl}|(H)$ for every $i$. For
every $i$ find a subgroup $E_i$ of $H$ such that $p_i^{\alpha_i}$
divides $|E_i|$. Let $E$ be a finite subgroup of $H$ containing
all the groups $E_i$. Clearly, $|F|$ divides $|E|$.
\end{pf}

\begin{Def}
A metric space is of {\it bounded geometry} if there is a number
$r>0$ and a function $c:\RR_+ \to \RR_+$ such that the
$r$-capacity (the maximal cardinality of $r$-discrete subset) of
every $\varepsilon$-ball does not exceed $c(\varepsilon)$.
\end{Def}

Notice that any countable group with proper left invariant metric
has bounded geometry.

A large scale analog $\MM^0$ of 0-dimensional Cantor set is
introduced in~\cite{Dran-Zar}: it is the set of all positive
integers with ternary expression containing 0's and 2's only (with
the metric from $\RR_+$): $ \MM^0=\{ \sum_{i=0}^\infty a_i3^i\mid
a_i=0,2\} $.

\begin{Prop}\cite[Theorem 3.11]{Dran-Zar}\label{M0 is universal}
The space $\MM^0$ is universal for proper metric spaces of
bounded geometry and of asymptotic dimension zero.
\end{Prop}

\begin{Prop}\label{M0 is equivalent to Z2}
The space $\MM^0$ is coarsely equivalent to
$\bigoplus\limits_{i=1}^\infty {\mathbb{Z}}_{2}$.
\end{Prop}

\begin{pf}
To define a map $f\colon \bigoplus\limits_{i=1}^\infty
{\mathbb{Z}}_{2}\to \MM^0$ we consider an element $p=p_1p_2\dots
p_n$ of the group $\bigoplus\limits_{i=1}^\infty {\mathbb{Z}}_{2}$ where
$p_j\in\{0,1\}={\mathbb{Z}}_{2}$ and put
$$f(p)=\sum_{i=1}^\infty 2p_i\cdot 3^{i-1}.$$
It is easy to check that the map $f$ is a coarse equivalence: for
any elements $p,q\in \bigoplus\limits_{i=1}^\infty {\mathbb{Z}}_{2}$ we
have
$$ 3^{d_{\mathcal{L}}(p,q)}\le
d_{\MM^0}(f(p),f(q))\le 3\cdot 3^{d_{\mathcal{L}}(p,q)}$$
\end{pf}

\begin{Rem}[cf. Proposition~\ref{ultrametric characterization}]
The proof above shows that the group $\bigoplus\limits_{i=1}^\infty
{\mathbb{Z}}_{2}$ with the ultrametric $3^{d_{\mathcal{L}}}$ is
bi-Lipschitz equivalent to the space $\MM^0$.
\end{Rem}

\begin{Prop}[cf.~{\cite[Theorem 4]{ProtasovMorphisms}}]
Let $G$ and $H$ be locally finite countable groups with proper
left invariant metrics. Then the metric space $G$ can be coarsely
embedded in the metric space $H$ (this map is not a homomorphism).
\end{Prop}

\begin{pf}
By Propositions~\ref{M0 is universal} and~\ref{M0 is equivalent to
Z2} the group $G$ can be coarsely embedded in the group $\oplus
{\mathbb{Z}}_{2}$. By Lemma~\ref{IsometricLocFinGroups} the group
$(\bigoplus {\mathbb{Z}}_{2}, d_{\mathcal{L}})$ embeds isometrically
into any group $(\bigoplus\limits_{i=1}^\infty {\mathbb{Z}}_{a_i},
d_{\mathcal{L}})$. Finally, the group $H$ is bi-uniformly
equivalent to a direct sum of cyclic groups by
Theorem~\ref{MetricOnLocallyFiniteGrps}.
\end{pf}

Let $G$ and $H$ be countable locally finite groups.
Using~\ref{MetricOnLocallyFiniteGrps} one can show that if
$$\sum_{p\text{-prime}}
\big| |\text{$p$\,-Syl}|(G)-|\text{$p$\,-Syl}|(H)\big|<\infty$$
then the groups $G$ and $H$ are coarsely equivalent. Is the
converse true?

\begin{Problem}\label{CharOfTorAb}
Classify countable abelian torsion groups up to coarse
equivalence.
\end{Problem}

Let us suggest a program to answer \ref{CharOfTorAb}. Notice that
any abelian torsion group is coarsely equivalent to a direct sum
of groups ${\mathbb{Z}}_p$ with $p$ being prime. Therefore the
following groups are of importance: ${\mathbb{Z}}^\infty_p$ (the
infinite direct sum of copies of ${\mathbb{Z}}_p$) and
$\bigoplus\limits_{p\in {\mathcal{P}}}{\mathbb{Z}}_{p^{n(p)}}$,
where $n(p)\ge 1$ for each $p\in {\mathcal{P}}$, ${\mathcal{P}}$
being a subset of primes.
\begin{Problem}
Suppose $\bigoplus\limits_{p\in {\mathcal{P}}}{\mathbb{Z}}_{p^{n(p)}}$ and
$\bigoplus\limits_{q\in {\mathcal{Q}}}{\mathbb{Z}}_{q^{m(q)}}$
are coarsely equivalent. Is the symmetric difference of ${\mathcal{P}}$ and ${\mathcal{Q}}$
finite? If so, does $n(p)$ equal $m(p)$ for all but finitely
many $p$?
\end{Problem}

\begin{Problem}
Suppose $\bigoplus\limits_{p\in {\mathcal{P}}}{\mathbb{Z}}^\infty_p$ and
$\bigoplus\limits_{q\in {\mathcal{Q}}}{\mathbb{Z}}^\infty_q$
are coarsely equivalent. Is ${\mathcal{P}}$ equal ${\mathcal{Q}}$?
\end{Problem}

Call two countable abelian torsion groups $G$ and $H$ {\it virtually isometric}
if there are subgroups of finite index $G'$ of $G$ and $H'$ of $H$
such that $G'$ is isometric to $H'$ for some choice of proper and invariant
metrics on $G'$ and $H'$. Notice virtually isometric groups are coarsely
equivalent.

\begin{Problem}
Suppose two countable abelian torsion groups $G$ and $H$
are coarsely equivalent. Are $G$ and $H$ virtually isometric?
\end{Problem}

\section{Examples of coarsely inequivalent ultrametric spaces}\label{section Examples}

In this section we construct uncountably many coarsely
inequivalent ultrametric spaces. Notice that any ultrametric space
has asymptotic dimension zero.

\begin{Def}
Let $(X,x_0)$ and $(Y,y_0)$ be pointed metric spaces. We define a
{\it metric wedge} $X\vee Y$ as the topological wedge of
these spaces with the following metric:
$$ d(z,z')\begin{cases}
d_X(z,z') \text{ if } z,z'\in X \\
d_Y(z,z') \text{ if } z,z'\in Y \\
\max\{d_X(z,x_0),d_Y(z',y_0)\} \text{ if $z\in X\setminus\{x_0\}$ and $z'\in Y \setminus\{y_0\}$}
\end{cases}
$$
\end{Def}

Similarly, one can define metric wedge of an arbitrary family of
pointed metric spaces (cf.~\cite[Example 2]{Bogatyi 1}
or~\cite[Theorem 2.2]{Bogatyi 2}).

The following Lemma is easy to prove.

\begin{Lem}
The metric wedge of any family of pointed ultrametric spaces is a
pointed ultrametric space.
\end{Lem}

If $X$ is a bounded ultrametric space of diameter less than $M$,
then {\it the cone} $\Cone(X,M)$ is obtained from $X$ by adding a
vertex $v$ and declaring $d(v,x)=M$ for all $x\in X$. $\Cone(X,M)$
is a pointed ultrametric space with the vertex $v$ being its base
point.

Our examples will be obtained by wedging cones over basic ultrametric spaces,
scaled copies of $0$-skeleta of simplices.

Given a set $\lambda$ of integers bigger than $1$, we create a list $X_i$, $i\ge 1$,
of spaces (called {\it islands}) satisfying the following conditions:
\begin{enumerate}
\item The cardinality $n_i$ of $X_i$ belongs to $\lambda$.
\item There is an integer $m_i\ge n_i$ such that $d(x,y)=m_i$
for all $x\ne y\in X_i$. Notice $m_i=\diam(X_i)$.
\item For each $m\ge n$ and $n\in \lambda$ the set
of islands $X_i$ such that $m=\diam(X_i)$ and $n=|X_i|$ is infinite.
\end{enumerate}

The wedge $X_\lambda$ of all $\Cone(X_i,k_i)$, where
$k_i=\sum\limits_{j\leq i}m_j$ (put $m_j=0$ for $j\leq 0$), is the
{\it $\lambda$-archipelago}. $k_i$ is the {\it separation} of
island $X_i$ in the $\lambda$-archipelago.

\begin{Prop}
If $\lambda_1\ne\lambda_2$, then the
$\lambda_1$-archipelago is not coarsely equivalent to the
$\lambda_2$-archipelago.
\end{Prop}
\begin{pf}

Let $X_1$ be a $\lambda_1$-archipelago, $X_2$ be a
$\lambda_2$-archipelago, and suppose that $f\colon X_{1}\to X_{2}$
and $g\colon X_{2}\to X_{1}$ are coarse equivalences such that the
maps $g\circ f$ and $f\circ g$ are $C$-close to the identity and
do not move the base points. Assume that the set
$\lambda_1\setminus\lambda_2$ is not empty and fix a number $n$ in
it.

There are three
parameters associated to an island in any archipelago: the size, the diameter,
and the separation. For simplicity, an $(n,N,S)$-island
contains $n$ points, is of diameter $N$, and separation $S$.
Notice $n\leq N\leq S$.

Let us explain the idea of the proof. Since the space $X_1$
contains a lot of $n$-point islands, we are going to choose an
$(n,N,S)$-island $P\subset X$ such that $f(P)$ is also
an $n$-point island in $X_2$. Since the archipelago $X_2$ has no
$n$-point islands, we get a contradiction. First we choose the size $N$ of the island $P$
to be so large that the map $f$ is injective on $P$ and the map
$g$ is injective on $f(P)$. Then we choose the separation $S$ of
the island $P$ to be so large that $f(P)$ is contained in some
island $Q$ in $X_2$ and $g(Q)$ is contained in some island in
$X_1$ (in fact, $g(Q)\subset P$).

Let us introduce some notations that we use in the rest of the
proof. Given a coarse equivalence $h\colon Y\to Z$ of metric
spaces we denote by $\rho_h$ and $\delta_h$ two real functions
such that $\rho_h(d_Y(y,y'))\le d_Z(h(y),h(y'))\le
\delta_h(d_Y(y,y'))$ for any $y,y'\in Y$. If one of the spaces
$Y$, $Z$ is unbounded then the other is also unbounded and
$\lim_{t\to\infty} \rho_h(t)=\infty=\lim_{t\to\infty}
\delta_h(t)$.

Fix an integer $N>C$ such that $\rho_f(N)>C$. Notice that since
$N>C$, any $(n,N,S)$-island $P\subset X_1$ is
$C$-discrete and $C$-separated from the rest of $X_{1}$. Therefore
the map $g\circ f$ is identity on $P$ and the map $f$ is injective
on $P$.

Clearly, the image $f(P)$ of any $(n,N,S)$-island
$P\subset X_1$ is $\delta_f(N)$-bounded in $X_2$ and therefore is
contained in one $\delta_f(N)$-component $Q$ of $X_2$. If the
island $P$ is $S$-separated in $X_1$, then its image $f(P)$ is at
least $\rho_f(S)$-far from the base point of $X_2$. We choose $S$
large enough to satisfy $\rho_f(S) > \delta_f(N)$ and thus to make
sure that the $\delta_f(N)$-component $Q$ containing $f(P)$ is an
island. Assume $Q$ is $(k,m,S')$-island where $m\le
\delta_f(N)$ and $k>n$ (recall that $f$ is injective on $P$).

Since $\rho_f(N)>C$, the image $f(P)$ is $C$-discrete and
therefore $m>C$. But then the map $f\circ g$ is identity on $Q$
and the map $g$ is injective on $Q$.

The image $g(Q)$ is $\delta_g(m)$-bounded and contains $P$. By
choosing $S$ to be greater than $\delta_g(\delta_f(N))$ we
guarantee that the island $P$ is more than $\delta_g(m)$-separated
from the rest of $X_{1}$, therefore the set $g(Q)$ is entirely in
$P$. Since $g$ is injective on $Q$, we must have $n\ge k$.
Contradiction.
\end{pf}

\begin{Cor}
There are uncountably many coarsely inequivalent asymptotically
0-dimensional subspaces of the ray $\RR_+$.
\end{Cor}
\begin{pf}
Due to Proposition~\ref{M0 is universal} it is sufficient to check
that every $\lambda$-archipelago $X$ is proper and has bounded geometry.

Given $R> 0$, a ball $\bar B(x,R)$ either coincides with $\bar B(x_0,R)$, where
$x_0$ is the center of the archipelago $X$, consists of $x$ only,
or is the island containing $x$ which has at most $R$ points in that case. Thus the number of
points in any ball $B(x,R)$ is bounded by some number depending on
$R$ only. This shows both $X$ being proper and of bounded
geometry.
\end{pf}


\begin{thebibliography}{99}

\bibitem{Banakh Protasov}
T.~Banakh, I.~Protasov, \emph{Ball stuctures and colorings of
graphs and groups}, Matem. Studii Monograph Series. 11, VNTL Publ.
2003, 148p.

\bibitem{Bogatyi 1}
Bogaty\u\i, S. A., \emph{A universal uniform rational ultrametric
on the space irrational numbers}, (Russian) Vestnik Moskov. Univ.
Ser. I Mat. Mekh. 2000, no. 6, 20--24; translation in Moscow Univ.
Math. Bull. {\bf 55} (2000), no. 6, 20--24 (2001).

\bibitem{Bogatyi 2}
Bogaty\u\i, S. A., \emph{Metrically homogeneous spaces}, (Russian)
Uspekhi Mat. Nauk {\bf 57} (2002), no. 2(344), 3--22; translation
in Russian Math. Surveys {\bf 57} (2002), no. 2, 221--240.

\bibitem{Brod-Dydak}
N.~Brodskiy, J.~Dydak, \emph{Coarse dimensions and partitions of
unity}, preprint math.GT/0506547.

\bibitem{Brod-Dydak-Higes-Mitra}
N.~Brodskiy, J.~Dydak, J.~Higes, A.~Mitra, \emph{Assouad-Nagata
dimension via Lipschitz extensions}, preprint math.MG/0601226.

\bibitem{Buyalo1}
S.~Buyalo, \emph{Asymptotic dimension of a hyperbolic space and
capacity dimension of its boundary at infinity}, Algebra i analis
(St. Petersburg Math. J.), {\bf 17} (2005), 70--95 (in Russian).

\bibitem{Buyalo Lebedeva}
S.~Buyalo, N.~Lebedeva, \emph{Dimension of locally and
asymptotically self-similar spaces}, preprint math.GT/0509433.

\bibitem{DavSem}
G. David, S. Semmes, \emph{Fractured fractals and broken dreams.
Self-similar geometry through metric and measure}. Oxford Lecture
Series in Mathematics and its Applications, 7. The Clarendon
Press, Oxford University Press, New York, 1997.

\bibitem{Dran asind}
A.~Dranishnikov, \emph{On asymptotic inductive dimension}, JP J.
Geom. Topol. {\bf 1} (2001), no.3, 239--247.

\bibitem{Dran asym top}
A.~Dranishnikov, \emph{Asymptotic topology}, Russian Math. Surveys
{\bf 55} (2000), no.6, 1085--1129.

\bibitem{Dran-Zar}
A.~Dranishnikov, M.~Zarichnyi, \emph{Universal spaces for
asymptotic dimension}, Topology Appl. {\bf 140} (2004), no.2-3,
203--225.

\bibitem{Gro asym invar}
M. Gromov, \emph{Asymptotic invariants for infinite groups}, in
Geometric Group Theory, vol. 2, 1--295, G.Niblo and M.Roller,
eds., Cambridge University Press, 1993.

\bibitem{Groot}
J. de~Groot, \emph{On a metric that characterizes dimension},
Canad. J. Math. {\bf 9} (1957), 511-–514.

\bibitem{Isbell book}
J.R.~Isbell, \emph{Uniform spaces}, Mathematical Surveys, No. 12
American Mathematical Society, Providence, R.I. 1964 xi+175 pp.

\bibitem{Lang-Sch Nagata dim}
U.~Lang, T.~Schlichenmaier, \emph{Nagata dimension, quasisymmetric
embeddings, and Lipschitz extensions}, IMRN International
Mathematics Research Notices (2005), no.58, 3625--3655.

\bibitem{Lemin ultrametrization}
A.~Lemin,  \emph{On ultrametrization of general metric spaces},
Proc. Amer. Math. Soc. {\bf 131} (2003), 979--989.

\bibitem{Lemin Lemin universal}
A.~Lemin, V.~Lemin, \emph{On a universal ultrametric space},
Topology Appl. {\bf 103} (2000), no.3, 339--345.

\bibitem{Nagata}
J. Nagata, \emph{On a relation between dimension and metrization},
Proc. Japan Acad. {\bf 32} (1956), 237-–240.

\bibitem{Ostrand}
P.~Ostrand, \emph{A conjecture of J. Nagata on dimension and
metrization}, Bull. Amer. Math. Soc. {\bf 71} (1965), 623--625.

\bibitem{ProtasovMorphisms}
I.V.~Protasov, \emph{Morphisms of ball's structures of groups and
graphs}, Ukrainian Math. J. {\bf 54} (2002), no. 6, 1027--1037.

\bibitem{ProtasovSurvey}
I.V.~Protasov, \emph{Survey of Balleans}, preprint.

\bibitem{Roe lectures}
J. Roe, \emph{Lectures on coarse geometry}, University Lecture
Series, 31. American Mathematical Society, Providence, RI, 2003.

\bibitem{Sha}
Y. Shalom, {\em Harmonic analysis, cohomology, and the large-scale geometry of amenable
groups}, Acta Math., 192 (2004), 119--185.

\bibitem{Smith}
 J. Smith, \emph{On Asymptotic Dimension of Countable Abelian Groups},
 preprint math.GR/0504447.

\bibitem{Urysohn}
P.S.~Urysohn, \emph{Sur un espace m\'{e}trique universel}, Bull.
Sci. Math. {\bf 51} (1927), 43--65 and 74--96.


\bibitem{MZar commun}
M.~Zarichnyi, private communication.

\end{thebibliography}
\end{document}